\documentclass[12pt]{amsart}
\usepackage[utf8]{inputenc}
\usepackage{xcolor}

\usepackage{subfiles}
\usepackage{hyperref}
\usepackage{cleveref}
\usepackage{amsmath}
\usepackage{multirow}
\usepackage{amssymb}
\usepackage{amsthm}
\usepackage{mathtools}
\usepackage{array}
\usepackage{thmtools}
\usepackage{thm-restate}
\usepackage[margin=1in]{geometry}
\usepackage{cleveref}
\usepackage{ytableau}
\usepackage{comment}
\usepackage{soul}
\usepackage{tikz}
\usetikzlibrary{positioning}

\usepackage{graphicx}
\usepackage{caption}
\usepackage{float}

\usepackage[maxnames=5]{biblatex}
\renewbibmacro{in:}{}
\DeclareFieldFormat[article]{citetitle}{#1}
\DeclareFieldFormat[article]{title}{#1}

\numberwithin{equation}{section}

\theoremstyle{definition}
\newtheorem{theorem}{Theorem}[section]
\newtheorem*{theorem*}{Theorem}

\newtheorem*{example*}{Example}
\newtheorem{lemma}[theorem]{Lemma}
\newtheorem*{lemma*}{Lemma}
\newtheorem{corollary}[theorem]{Corollary}
\newtheorem*{corollary*}{Corollary}

\newtheorem*{definition*}{Definition}
\newtheorem{proposition}[theorem]{Proposition}
\newtheorem*{proposition*}{Proposition}

\newtheorem*{remark*}{Remark}

\usepackage{ mathrsfs }

\title[Set partitions that require maximum sorts through the $aba-$avoiding stack]{On the set partitions that require \\ maximum sorts through the $aba-$avoiding stack}

\author{Yunseo Choi}\address{\textsc{Y. Choi}, Harvard University,
    Cambridge, MA, 02138} \email{ychoi@college.harvard.edu}

\author{Katelyn Gan}\address{\textsc{K. Gan}, Sage Hill School,
    Newport Beach, CA, 92657} \email{katelyngan77@gmail.com}

\author{Andrew Li}\address{\textsc{A. Li}, Highland Park High School, Dallas, TX, 75205} \email{andrewli10062006@gmail.com}

\author{Tiffany Zhu}\address{\textsc{T. Zhu}, The Harker School, San Jose, CA, 95129} \email{26tiffanyz@gmail.com}

\begin{document}

\maketitle

\begin{abstract} 
Recently, Xia introduced a deterministic variation $\phi_{\sigma}$ of Defant and Kravitz's stack-sorting maps for set partitions and showed that any set partition $p$ is sorted by $\phi^{N(p)}_{aba}$, where $N(p)$ is the number of distinct alphabets in $p$. Xia then asked which set partitions $p$ are not sorted by $\phi_{aba}^{N(p)-1}$. In this note, we prove that the minimal length of a set partition $p$ that is not sorted by $\phi_{aba}^{N(p)-1}$ is $2N(p)$. Then we show that there is only one set partition of length $2N(p)$ and ${{N(p) + 1} \choose 2} + 2{N(p) \choose 2}$  set partitions of length $2N(p)+1$ that are not sorted by $\phi_{aba}^{N(p)-1}$.
\end{abstract}

\section{Introduction} 
\label{intro}
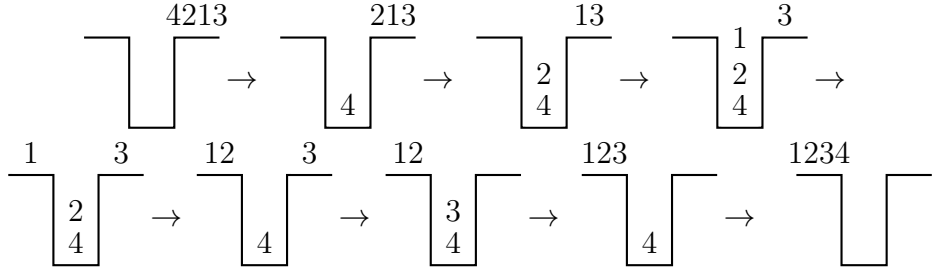
\begin{figure}[h]\begin{center}
\begin{tikzpicture}[scale=0.6]
				\draw[thick] (0,0) -- (1,0) -- (1,-2) -- (2,-2) -- (2,0) -- (3,0);
				\node[fill = white, draw = white] at (2.5,.5) {4213};
				\node[fill = white, draw = white] at (1.5,-1.5) {};
				\node[fill = white, draw = white] at (1.5,-.8) {};
				\node[fill = white, draw = white] at (3.5,-1) {$\rightarrow$};
              \end{tikzpicture}
              \begin{tikzpicture}[scale=0.6]
				\draw[thick] (0,0) -- (1,0) -- (1,-2) -- (2,-2) -- (2,0) -- (3,0);
				\node[fill = white, draw = white] at (2.5,.5) {213};
				\node[fill = white, draw = white] at (1.5,-1.5) {4};
                \node[fill = white, draw = white] at (1.5,-.8) {};
				\node[fill = white, draw = white] at (.5,.5) {};
             \node[fill = white, draw = white] at (3.5,-1) {$\rightarrow$}; \end{tikzpicture}
         \begin{tikzpicture}[scale=0.6]
				\draw[thick] (0,0) -- (1,0) -- (1,-2) -- (2,-2) -- (2,0) -- (3,0);
				\node[fill = white, draw = white] at (2.5,.5) {13};
				\node[fill = white, draw = white] at (1.5,-1.5) {4};
				\node[fill = white, draw = white] at (1.5,-.8) {2};
				\node[fill = white, draw = white] at (3.5,-1) {$\rightarrow$};
              \end{tikzpicture} 
              \begin{tikzpicture}[scale=0.6]
				\draw[thick] (0,0) -- (1,0) -- (1,-2) -- (2,-2) -- (2,0) -- (3,0);
				\node[fill = white, draw = white] at (2.5,.5) {3};
				\node[fill = white, draw = white] at (1.5,-1.5) {4};
				\node[fill = white, draw = white] at (1.5,-.8) {2};
                \node[fill = white, draw = white] at (1.5,0) {1};
             	\node[fill = white, draw = white] at (3.5,-1) {$\rightarrow$}; \end{tikzpicture}
              \begin{tikzpicture}[scale=0.6]
				\draw[thick] (0,0) -- (1,0) -- (1,-2) -- (2,-2) -- (2,0) -- (3,0);
				\node[fill = white, draw = white] at (2.5,.5) {3};
				\node[fill = white, draw = white] at (1.5,-1.5) {4};
                 \node[fill = white, draw = white] at (1.5,-.8) {2};
				\node[fill = white, draw = white] at (.5,.5) {1};
             	\node[fill = white, draw = white] at (3.5,-1) {$\rightarrow$}; \end{tikzpicture}\begin{tikzpicture}[scale=0.6]
				\draw[thick] (0,0) -- (1,0) -- (1,-2) -- (2,-2) -- (2,0) -- (3,0);
				\node[fill = white, draw = white] at (2.5,.5) {3};
				\node[fill = white, draw = white] at (1.5,-1.5) {4};
                 \node[fill = white, draw = white] at (1.5,-.8) {};
                 \node[fill = white, draw = white] at (1.5,-.1) {};
				\node[fill = white, draw = white] at (.5,.5) {12};
             	\node[fill = white, draw = white] at (3.5,-1) {$\rightarrow$}; \end{tikzpicture}\begin{tikzpicture}[scale=0.6]
				\draw[thick] (0,0) -- (1,0) -- (1,-2) -- (2,-2) -- (2,0) -- (3,0);
				\node[fill = white, draw = white] at (2.5,.5) {};
				\node[fill = white, draw = white] at (1.5,-1.5) {4};
                 \node[fill = white, draw = white] at (1.5,-.8) {3};
				\node[fill = white, draw = white] at (.5,.5) {12};
             	\node[fill = white, draw = white] at (3.5,-1) {$\rightarrow$}; \end{tikzpicture}\begin{tikzpicture}[scale=0.6]
				\draw[thick] (0,0) -- (1,0) -- (1,-2) -- (2,-2) -- (2,0) -- (3,0);
				\node[fill = white, draw = white] at (2.5,.5) {};
				\node[fill = white, draw = white] at (1.5,-1.5) {4};
                 \node[fill = white, draw = white] at (1.5,-.8) {};
				\node[fill = white, draw = white] at (.5,.5) {123};
             	\node[fill = white, draw = white] at (3.5,-1) {$\rightarrow$}; \end{tikzpicture}
              \begin{tikzpicture}[scale=0.6]
				\draw[thick] (0,0) -- (1,0) -- (1,-2) -- (2,-2) -- (2,0) -- (3,0);
				\node[fill = white, draw = white] at (.5,.5) {1234};\end{tikzpicture}
\end{center}
\caption{West's stack-sorting map $s$ on $\pi = 4213$}
\label{Westmap}
\end{figure}

In 1973, Knuth \cite{knuth1997art} introduced a non-deterministic stack-sorting machine that at each step, either pushes the leftmost remaining entry of the input permutation into the stack or pops the topmost entry of the stack. In 1990, West \cite{WEST1990} modified Knuth's stack-sorting machine so that it is deterministic. In West's deterministic stack-sorting map $s$, the input permutation is sent through a stack in a right-greedy manner, while insisting that the stack is increasing when read from top to bottom (see for example, \Cref{Westmap}). It is a classical result that $s^{n-1}(\pi) = \mathrm{id}$ for any $\pi \in S_n$.

In 2020, Cerbai, Claesson, and Ferrari \cite{CERBAI2020105230} extended West's stack-sorting map $s$ to $s \circ s_{\sigma}$, where the map $s_{\sigma}$ sends the input permutation through a stack in a right greedy manner, while maintaining that the stack avoids subsequences that are order-isomorphic to some permutation $\sigma$ (Note that $s_{21} = s$). In the following year, Berlow \cite{BERLOW2021112571} generalized $s_{\sigma}$ to $s_{T}$, in which the stack must simultaneously avoid subsequences that are order isomorphic to any of the permutations in the set $T$, while Defant and Zheng \cite{DEFANT2021102192} generalized $s_{\sigma}$ to $s_{\overline{\sigma}}$, in which the stack must avoid substrings that are order isomorphic to $\sigma$ at all times.  

\begin{figure}[h]\begin{center}
\scalebox{.9}{
\begin{tikzpicture}[scale=0.6]
				\draw[thick] (-0.5,0) -- (1,0) -- (1,-2) -- (2,-2) -- (2,0) -- (3.5,0);
				\node[fill = white, draw = white] at (2.75,.5) {abcac};
				\node[fill = white, draw = white] at (1.5,-1.5) {};
				\node[fill = white, draw = white] at (1.5,-.8) {};
				\node[fill = white, draw = white] at (4,-1) {$\longrightarrow$};
              \end{tikzpicture}} \hspace{-5mm}
              \scalebox{.9}{
\begin{tikzpicture}[scale=0.6]
				\draw[thick] (-0.5,0) -- (1,0) -- (1,-2) -- (2,-2) -- (2,0) -- (3.5,0);
				\node[fill = white, draw = white] at (2.75,.5) {bcac};
				\node[fill = white, draw = white] at (1.5,-1.5) {a};
				\node[fill = white, draw = white] at (1.5,-.8) {};
				\node[fill = white, draw = white] at (4,-1) {$\longrightarrow$};
              \end{tikzpicture}} \hspace{-5mm}
              \scalebox{.9}{
\begin{tikzpicture}[scale=0.6]
				\draw[thick] (-0.5,0) -- (1,0) -- (1,-2) -- (2,-2) -- (2,0) -- (3.5,0);
				\node[fill = white, draw = white] at (2.75,.5) {cac};
				\node[fill = white, draw = white] at (1.5,-1.5) {a};
				\node[fill = white, draw = white] at (1.5,-.8) {b};
				\node[fill = white, draw = white] at (4,-1) {$\longrightarrow$};
              \end{tikzpicture}} \hspace{-5mm}
              \scalebox{.9}{
\begin{tikzpicture}[scale=0.6]
				\draw[thick] (-0.5,0) -- (1,0) -- (1,-2) -- (2,-2) -- (2,0) -- (3.5,0);
				\node[fill = white, draw = white] at (2.75,.5) {ac};
				\node[fill = white, draw = white] at (1.5,-1.5) {a};
				\node[fill = white, draw = white] at (1.5,-.8) {b};
                     \node[fill = white, draw = white] at (1.5,-.1) {c};
				\node[fill = white, draw = white] at (4,-1) {$\longrightarrow$};
              \end{tikzpicture}} \hspace{-5mm}
              \scalebox{.9}{
\begin{tikzpicture}[scale=0.6]
				\draw[thick] (-0.5,0) -- (1,0) -- (1,-2) -- (2,-2) -- (2,0) -- (3.5,0);
				\node[fill = white, draw = white] at (2.75,.5) {ac};
				\node[fill = white, draw = white] at (1.5,-1.5) {a};
				\node[fill = white, draw = white] at (1.5,-.8) {b};
				\node[fill = white, draw = white] at (4,-1) {$\longrightarrow$};
    \node[fill = white, draw = white] at (.25,.5) {c};
              \end{tikzpicture}} \hspace{-5mm}
              \scalebox{.9}{
\begin{tikzpicture}[scale=0.6]
				\draw[thick] (-0.5,0) -- (1,0) -- (1,-2) -- (2,-2) -- (2,0) -- (3.5,0);
				\node[fill = white, draw = white] at (2.75,.5) {ac};
				\node[fill = white, draw = white] at (1.5,-1.5) {a};
				\node[fill = white, draw = white] at (4,-1) {$\longrightarrow$};
    \node[fill = white, draw = white] at (.25,.5) {cb};
              \end{tikzpicture}} \hspace{-5mm}
                            \scalebox{.9}{
\begin{tikzpicture}[scale=0.6]
				\draw[thick] (-0.5,0) -- (1,0) -- (1,-2) -- (2,-2) -- (2,0) -- (3.5,0);
				\node[fill = white, draw = white] at (2.75,.5) {c};
				\node[fill = white, draw = white] at (1.5,-1.5) {a};
    \node[fill = white, draw = white] at (1.5,-.8) {a};
				\node[fill = white, draw = white] at (4,-1) {$\longrightarrow$};
    \node[fill = white, draw = white] at (.25,.5) {cb};
              \end{tikzpicture}} \hspace{-5mm}
\scalebox{.9}{
\begin{tikzpicture}[scale=0.6]
				\draw[thick] (-0.5,0) -- (1,0) -- (1,-2) -- (2,-2) -- (2,0) -- (3.5,0);
    				\node[fill = white, draw = white] at (1.5,-.1) {c};
				\node[fill = white, draw = white] at (1.5,-1.5) {a};
    \node[fill = white, draw = white] at (1.5,-.8) {a};
				\node[fill = white, draw = white] at (4,-1) {$\longrightarrow$};
    \node[fill = white, draw = white] at (.25,.5) {cb};
              \end{tikzpicture}} \hspace{-5mm}
                                          \scalebox{.9}{
\begin{tikzpicture}[scale=0.6]
				\draw[thick] (-0.5,0) -- (1,0) -- (1,-2) -- (2,-2) -- (2,0) -- (3.5,0);
    				\node[fill = white, draw = white] at (1.5,-0.8) {a};
				\node[fill = white, draw = white] at (1.5,-1.5) {a};
				\node[fill = white, draw = white] at (4,-1) {$\longrightarrow$};
    \node[fill = white, draw = white] at (.25,.5) {cbc};
              \end{tikzpicture}} \hspace{-5mm}
                                          \scalebox{.9}{
\begin{tikzpicture}[scale=0.6]
				\draw[thick] (-0.5,0) -- (1,0) -- (1,-2) -- (2,-2) -- (2,0) -- (3.5,0);
				\node[fill = white, draw = white] at (1.5,-1.5) {a};
				\node[fill = white, draw = white] at (4,-1) {$\longrightarrow$};
    \node[fill = white, draw = white] at (.25,.5) {cbca};
              \end{tikzpicture}} \hspace{-5mm}
            \scalebox{.9}{
\begin{tikzpicture}[scale=0.6]
				\draw[thick] (-0.5,0) -- (1,0) -- (1,-2) -- (2,-2) -- (2,0) -- (3.5,0);
    \node[fill = white, draw = white] at (.25,.5) {cbcaa};
              \end{tikzpicture}}
\end{center}
\caption{Xia's stack-sorting map $\phi_{aba}$ on $p= abcac$}
\label{xiamap}
\end{figure}
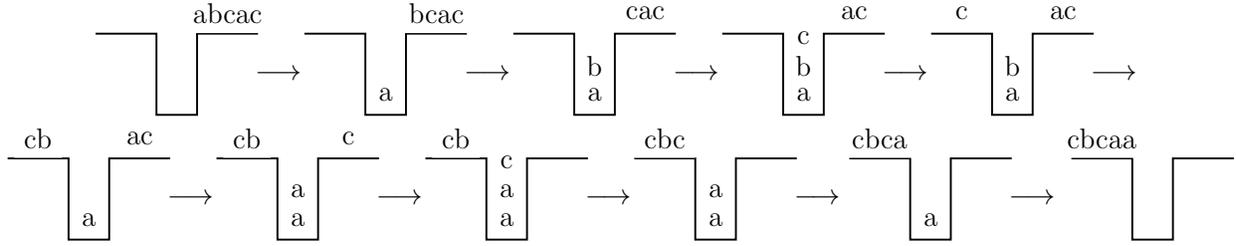

More recently, in 2022, Defant and Kravitz \cite{defant2022footsorting} generalized Knuth's non-deterministic stack-sorting-machine \cite{knuth1997art} to \textit{set partitions}, which are sequences of (possibly repeated) letters from some set of alphabets $A$. In 2023, Xia \cite{xia2023deterministic} introduced deterministic variations, $\phi_{\sigma}$ and $\phi_{\overline{\sigma}}$, of Defant and Kravitz's stack-sorting map for set partitions \cite{defant2022footsorting} as did West \cite{WEST1990} of Knuth's stack-sorting machine \cite{knuth1997art}. A set partition is said to be \textit{sorted} if all occurrences of the same letter appear consecutively on the set partition, and two set partitions $p = p_1 p_2 \cdots p_n$ and $q = q_1 q_2 \cdots q_n$ are \textit{equivalent} if there exists some bijection $f: A \to A$ such that $q = f(p_1) f(p_2) \cdots f(p_n)$. In Xia's deterministic stack-sorting map $\phi_{\sigma}$ for set partitions, the input set partition is sent through a stack in a right-greedy manner, while insisting that the stack avoids subsequences that are equivalent to the set partition $\sigma$ (see for example, \Cref{xiamap}). 

In addition to introducing $\phi_{\sigma}$, Xia \cite[Proposition 5.2]{xia2023deterministic} showed that $\phi_{aba}$ is the only $\phi_{\sigma}$ that eventually sorts all set partitions. Then in \cite[Theorem 3.2]{xia2023deterministic}, Xia showed that any set partition $p$ is sorted after applying $\phi^{N(p)}_{aba}$ and demonstrated the sharpness of her bound by proving that $p=(a_1 a_2 \cdots a_{N(p)})^2$ is not sorted after applying $\phi^{N(p)-1}_{aba}$ for any $N(p) \geq 3$. Finally, Xia asked \cite[Question 1]{xia2023deterministic} about which set partitions $p$ are not sorted after applying $\phi^{N(p)-1}_{aba}$. We first answer Xia's question with the restriction that $|p| \leq 2N(p)$. 

\begin{theorem} \label{main}
    If set partition $p$ satisfies $|p| \leq 2N(p)$ for some $N(p) \geq 3$ and is not sorted after applying $\phi^{N(p)-1}_{aba}$, then $p$ is equivalent to $(a_1 a_2 \cdots a_{N(p)})^2.$
\end{theorem}

\Cref{main} proves that for any fixed $N(p) \geq 3$, the example that Xia \cite{xia2023deterministic} used in her Theorem 3.2 is, up to equivalence, the only shortest set partition $p$ that is not sorted after applying $\phi^{N(p)-1}_{aba}$. In \Cref{main2}, we enumerate the next-minimal set partitions $p$ that are not sorted by $\phi_{aba}^{N(p)-1}(p)$---the set partitions of length $2N(p)+1$ that are not sorted after applying $\phi_{aba}^{N(p)-1}(p)$. 

\begin{theorem} \label{main2}
For a fixed $N(p) \geq 3$, the number of inequivalent set partitions $p$ that satisfy $|p| = 2N(p) + 1$ and are not sorted after applying $\phi_{aba}^{N(p)-1}(p)$ is ${N(p)+1 \choose 2} +  2 {N(p) \choose 2}$. 
\end{theorem}

 The rest of this note is organized as follows. In \Cref{prelim}, we establish the preliminaries. In \Cref{proofs}, we prove \Cref{main} and \Cref{main2}. 

\section{Preliminaries}
\label{prelim}
Let $A$ be an infinite set of alphabets. In this note, we use $a_1, a_2, a_3, \ldots$ or the standard Latin alphabets $a, b, c, \ldots$ to refer to alphabets of $A$. Unless otherwise specified, $a_1, a_2, a_3, \ldots$ are distinct alphabets of $A$. First, for (possibly empty) set partitions $s_1, s_2, \ldots, s_n$, let $s_1 s_2 \cdots s_n$ be their concatenation. Next, for a (possibly empty) set partition $p$, let $|p|$ be its length, and let $p^{m} = \underbrace{p p \cdots p}_{m \text{ times}}$. In addition, for a (possibly empty) set partition $p=p_1p_2 \cdots p_{|p|}$, let $p_{[i:j]} = p_i p_{i+1} \cdots p_j$. Now, let the \textit{reverse} of a set partition $p$ be $\mathrm{r}(p)=p_{|p|} p_{|p|-1} \cdots p_1$. So, for example, if $p = abcac$, then $\mathrm{r}(p) = cacba$. Furthermore, for $p=p_1p_2 \cdots p_{|p|}$ and $a \in A$, say that $a \in p$ if there exists some $i$ such that $p_i=a$, and for a set of alphabets $B \subseteq A$, let $\mathrm{ind}_{p}(B)$ be the set of $i$ such that $p_i \in B$. If $|B|=1$, then we omit the brackets around the set $B$ when writing $\mathrm{ind}_{p}(B)$. For example, if $p=a_1 a_2 a_2 a_3 a_1 a_1$, then $\mathrm{ind}_{p}(a_1) = \{1,5,6\}$, and $\mathrm{ind}_{p}(\{a_1, a_3\}) = \{1,4,5,6\}$. Let the $i^{\text{th}}$ smallest number in the set $\mathrm{ind}_{p}(B)$ be $\mathrm{ind}_{p}^{i}(B)$. Next, for any $p$, let $\mathrm{mcount}(p) = \max_{a \in A} |\mathrm{ind}_{p}(a)|$. For example, $\mathrm{mcount}(p) = 2$ for $p = a_1a_2a_3a_1a_3$. Now, for $\{a_j, a_k\} \subseteq p$ that satisfy $|\mathrm{ind}_{p}(a_j)| =|\mathrm{ind}_{p}(a_k)| =2 $, say that $a_j$ and $a_j$ are \textit{crossing} in $p$ if either $\mathrm{ind}_{p}(a_j)$ or $\mathrm{ind}_{p}(a_k)$ is $\{\mathrm{ind}_{p}^{1}(\{a_j, a_k\}), \mathrm{ind}_{p}^{3}(\{a_j, a_k\})\}$. For example, $a_1$ and $a_2$ are crossing in $p =a_1a_2a_1a_3a_2a_3$ but $a_1$ and $a_3$ are not. Also, as in Xia \cite{xia2023deterministic}, say that an alphabet in $p$ is \textit{clumped} in $p$ if all instances of the alphabet appear consecutively in $p$. 
 Let $C(p)$ be the number of clumped alphabets in $p$, and let $\mathrm{nc}(p)$ be the leftmost alphabet in $p$ that is not clumped in $p$. For example, in $p =a_1 a_1 a_1 a_2 a_3 a_4 a_2 a_4$, the alphabets $a_1$ and $a_3$ are clumped, so $C(p) = 2$ and $\mathrm{nc}(p) = a_2$. Note that $p$ is sorted if and only if $C(p) = N(p)$. Now, every set partition $p$ can be uniquely written as $p = a_1^{\ell_1} a_2^{\ell_2} \cdots a_{m}^{\ell_m}$ for some possibly repeating set of alphabets $a_1, a_2, \ldots, a_m$ such that $a_{i} \ne a_{i+1}$ for all $1 \leq i \leq m-1$ and $\ell_i > 0$ for all $1 \leq i \leq m$. Then let the \textit{truncation} of a set partition $p$ be $\mathrm{trunc}(p) = a_1 a_2 \cdots b_{m}$. For example, if $p = a_1 a_1 a_1 a_2a_2 a_1 a_1 a_3$, then $\mathrm{trunc}(p) = a_1 a_2 a_1 a_3$. We end this section by citing a lemma and a corollary in Xia \cite{xia2023deterministic}. 
\begin{lemma}[Xia {\cite[Lemma 3.1]{xia2023deterministic}}] \label{lem1}
    Let $p = p_1^{\ell_1} s_1 p_1^{\ell_2} \cdots p_1^{\ell_{m}} s_{m} p_1^{\ell_{m+1}}$ for $\ell_1, \ell_2, \ldots, \ell_{m} > 0$ and $\ell_{m+1} \geq 0$ such that $p_1 \not\in s_i$ for all $1 \leq i \leq m$. Then
    \begin{align*}
        \phi_{aba}(p) = \phi_{aba}(s_1) \phi_{aba}(s_2) \cdots \phi_{aba}(s_m) p_1^{\ell_1 + \ell_2 + \cdots + \ell_{m+1}}. 
    \end{align*}
\end{lemma}

Now, it follows as a corollary of \Cref{lem1} that if $p$ is not sorted, then $C(\phi_{aba}(p)) > C(p)$, because $\mathrm{nc}(p)$ is not clumped in $p$ but is clumped in $\phi_{aba}(p)$. 

\begin{corollary}[Xia {\cite[Proof of Theorem 3.2]{xia2023deterministic}}]\label{cor1}
    If $p$ is not sorted, then $C(\phi_{aba}(p)) > C(p)$. 
\end{corollary}

The following corollary then follows immediately from \Cref{cor1}. 

\begin{corollary}\label{cor2}
    If $p$ is not sorted by $\phi_{aba}^{N(p)-1}$, then $C(\phi_{aba}^{i}(p)) = i$ for all $0 \leq i \leq N(p)$. 
\end{corollary}

\section{Proofs of the Main Results}
\label{proofs}
To prove \Cref{main}, we first note that the following proposition follows directly from the definition of truncation. 
\begin{proposition} \label{trunc}
    For any $p$, it holds that $\mathrm{trunc}(\phi_{aba}(p))=\mathrm{trunc}(\phi_{aba}(\mathrm{trunc}(p)))$.
\end{proposition}

We now prove  \Cref{main} through \Cref{lem1}, \Cref{cor2}, and \Cref{trunc}. 

\begin{proof}[Proof of \Cref{main}]
    By Xia \cite[Theorem 3.2]{xia2023deterministic}, any set partition that is equivalent to $(a_1 a_2 \cdots a_{N(p)})^2$ is not sorted after applying $\phi^{N(p)-1}_{aba}$ for $N(p) \geq 3$. It thus suffices to show that if $p$ satisfies $|p| \leq 2N(p)$ and is not sorted after applying $\phi^{N(p)-1}_{aba}$, then it is equivalent to $(a_1 a_2 \cdots a_{N(p)})^2,$ towards which, we induct on $N(p)$. 
    
    The statement clearly holds for $N(p) =3$. Now, suppose that $N(p) \geq 4$ and that if $p$ satisfies $|p| \leq 2N(p)-2$ and is not sorted after applying $\phi^{N(p)-2}_{aba}$, then it is equivalent to $(a_1 a_2 \cdots a_{N(p)-1})^2$. First, by \Cref{cor2}, $C(\phi_{aba}^{0}(p)) = C(p) = 0$, and so every $a \in p$ must satisfy $|\mathrm{ind}_{p}(a)| \geq 2$. But because $|p|\leq 2N(p)$, it must be that $|\mathrm{ind}_{p}(a)| = 2$ for all $a \in p$. Now, let $p= p_1 s_1 p_1 s_2$ for some set partitions $s_1$ and $s_2$. Then because $C(\phi_{aba}(p))=1$, each $a (\ne p_1) \in p$ satisfies $a \in s_1$ and $a \in s_2$; otherwise, by \Cref{lem1}, at least one of $\mathrm{nc}(s_1)$ or $\mathrm{nc}(s_2)$ are clumped in $\phi_{aba}(p)$ in addition to $p_1$, which negates \Cref{cor2} for $i=1$.

     Now, because all $a (\ne p_1) \in p$ satisfy $a \in s_1$ and $a \in s_2$, if $p$ is not equivalent to $(a_1 a_2 \cdots a_{N(p)})^2$, then some $a_j$ and $a_k$ must not be crossing in $p$. Furthermore, by \Cref{lem1}, the same $a_j$ and $a_k$ must not be crossing in $\phi_{aba}(p)$ as well. Now, by \Cref{trunc}, the set partition $q= \phi_{aba}(p)_{[1: 2N(p)-2]}$ satisfies $|q| = 2N(p)-2 = 2N(q)$, and $\phi_{aba}^{N(q)-1}(q)$ must not be sorted; otherwise, $\phi_{aba}^{N(p)-1}(p)$ will be sorted. Thus, by the induction hypothesis, $a_j$ and $a_k$ must be crossing in $q= \phi_{aba}(p)_{[1: 2N(p)-2]}$. Therefore, $p$ must be equivalent to $(a_1 a_2 \cdots a_{N(p)})^2$. 
\end{proof}

Next, we prove auxiliary lemmas that lead up to \Cref{main2}. First, for $|p|=2N(p)+1$ and $p$ that is not sorted after applying $\phi_{aba}^{N(p)-1}$, we prove that $\mathrm{ind}_{p}(a)=2$ for all but one $a \in p$ and $\mathrm{ind}_{p}(a_{*})=3$ for exactly one $a_{*} \in p$.

\begin{lemma} \label{char}
    If $p$ satisfies $|p| = 2N(p)+1$ and is not sorted after applying $\phi_{aba}^{N(p)-1}$, then there exists exactly one $a_{*} \in p$ such that $|\mathrm{ind}_{p}(a_{*})| = 3$, and for any other $a \in p$, it holds that $|\mathrm{ind}_{p}(a)| = 2$. 
\end{lemma}
\begin{proof}
    By \Cref{cor2}, $C(p) = 0$. Thus, $|\mathrm{ind}_{p}(a)| \geq 2$ for all $a\in p$. Then because $|p| = 2N(p)+1$, all but one $a \in p$ must satisfy $|\mathrm{ind}_{p}(a)| = 2$ and one $a_{*} \in p$ must satisfy $|\mathrm{ind}_{p}(a_{*})| = 3$.
\end{proof}

Next, we show that if a set partition $p$ satisfies the statement of \Cref{main2} and in addition $|\mathrm{ind}_{p}(p_1)|=2$, then the subsequence $p_{[1:\mathrm{ind}_p^2(p_1)]}$ or $p_{[\mathrm{ind}_p^2(p_1): |p|]}$ contains at most 2 of each alphabet.  
\begin{lemma} \label{theorem2lem1}
    If $p$ satisfies $|p| = 2N(p)+1$, is not sorted after applying $\phi_{aba}^{N(p)-1}$, and satisfies $|\mathrm{ind}_{p}(p_1)|=2$, then $\mathrm{mcount}(p_{[1:\mathrm{ind}_p^2(p_1)]})\leq 2$. Similarly,  $\mathrm{mcount}(p_{[\mathrm{ind}_p^2(p_1): |p|]}) \leq 2$.
\end{lemma}
\begin{proof}
Let $p = p_1 s_1 p_1 s_2$, and let $a_{*} \in p$ be as in the statement of \Cref{char}. By \Cref{char}, it suffices to show that $|\mathrm{ind}_{s_1}(a_{*})| < 3$ and that $|\mathrm{ind}_{s_2}(a_{*})| < 3$. Now, if $|\mathrm{ind}_{s_1}(a_{*})| = 3$, then $\mathrm{nc}(s_1)$ gets clumped in $\phi_{aba}(p)$. Thus, $C(\phi_{aba}(p)) \geq 2$, a contradiction to \Cref{cor2} for $i=1$. Thus, $|\mathrm{ind}_{s_1}(a_{*})| < 3$. Similarly,  $|\mathrm{ind}_{s_2}(a_{*})| < 3$. 
\end{proof}

Next, we count the inequivalent set partitions $p$ that satisfy the statement of \Cref{main2} and in addition, $\mathrm{mcount}(p_{[\mathrm{ind}_p^2(p_1)+1: |p|]}) =2$. 

\begin{lemma} \label{theorem2lem2}
The number of inequivalent $p$ that satisfies $|p| = 2N(p)+1$, is not sorted after applying $\phi_{aba}^{N(p)-1}$, and satisfies $|\mathrm{ind}(p_1)|=\mathrm{mcount}(p_{[\mathrm{ind}_p^2(p_1)+1: |p|]}) =2$ is ${N(p) \choose 2}$. 
\end{lemma}

\begin{proof}
     Let $a_{*}$ be as in the statement of \Cref{char}, and let  $p=p_1s_1p_1s_2a_{*}s_3a_{*}s_4$ for some set partitions $s_1, s_2, s_3,$ and $s_4$. In additon, let $S = \{s_1, s_2, s_3, s_4\}$. Now, for any $s \in S$, if some $a \in s$ satisfies $|\mathrm{ind}_{s}(a)| = 2$, then $\mathrm{nc}(s)$ is clumped in $\phi_{aba}(p)$ by \Cref{lem1}. But this negates \Cref{cor2} for $i=1$. Therefore, all $a \in s$ for each $s \in S$ must satisfy $|\mathrm{ind}_{s}(a)| = 1$. Furthermore, if some $a \in s_2$ also satisfies $a \in s_3$ (resp. $s_4$), then $\mathrm{nc}(s_2 a_{*} s_3 a_{*} s_4)$ is clumped in $\phi_{aba}(p)$. But if so, then $C(\phi_{aba}(p)) > 1$, which negates \Cref{cor2} for $i=1$; thus, $s_2 \cap s_3 = \emptyset$ (resp, $s_2 \cap s_4 = \emptyset$). As a result, $\phi_{aba}(p) = \mathrm{r}(s_1) \mathrm{r}(s_3) \mathrm{r}(s_4) {a_{*}}^2 \mathrm{r}(s_2) p_1^2$ and so, $\mathrm{trunc}(\phi_{aba}(p)) = \mathrm{r}(s_1) \mathrm{r}(s_3) \mathrm{r}(s_4) a_{*} \mathrm{r}(s_2) p_1$.
    
    Next, by \Cref{trunc}, if $p$ is not sorted by $\phi_{aba}^{N(p)-1},$ then $\mathrm{r}(s_1) \mathrm{r}(s_3) \mathrm{r}(s_4) a_{*} \mathrm{r}(s_2)$ must not be sorted by $\phi_{aba}^{N(p)-2}$. Thus, by \Cref{main}, it must be that  $$\mathrm{r}(s_1) \mathrm{r}(s_3) \mathrm{r}(s_4) a_{*} \mathrm{r}(s_2) = (\phi_{aba}(p)_1 \phi_{aba}(p)_2 \cdots \phi_{aba}(p)_{N(p)-1})^2.$$

     Now, because $a_{*} \in s_1$ by \Cref{theorem2lem1} and $s_2 \cap (s_3 \cup s_4) = \emptyset$, it must be that $|r(s_1)| \geq N(p)-1$. But because each $a \in s_1$ satisfies $|\mathrm{ind}_{a}(s_1)|=1$, it holds that $|r(s_1)| \leq N(p)-1$. Thus, $|r(s_1)| = N(p)-1$. As a result, $$\mathrm{r}(s_3) \mathrm{r}(s_4) a_{*} \mathrm{r}(s_2) = \mathrm{r}(s_1) = {\mathrm{r}(p_2 p_3 \cdots p_{N(p)})}.$$ Thus, fixing $|s_2|,$ $|s_3|,$ and $|s_4|$ such that $|s_2| +|s_3| + |s_4| = N(p)-1$ fixes $p$.
\end{proof}

Similarly, we count the inequivalent set partitions $p$ that satisfy the statement of \Cref{main2} and in addition, $\mathrm{mcount}(p_{[1:\mathrm{ind}_p^2(p_1)-1]}) =2$.

\begin{lemma}\label{theorem2lem3}
    The number of inequivalent $p$ that satisfies $|p| = 2N(p)+1$, is not sorted after applying $\phi_{aba}^{N(p)-1}$, and satisfies $|\mathrm{ind}(p_1)|=\mathrm{mcount}(p_{[1:\mathrm{ind}_p^2(p_1)-1]}) =2$ is ${N(p) \choose 2}$. 
\end{lemma}

\begin{proof}
    Let $a_{*}$ be as in the statement of \Cref{char}, and let $p = p_1 s_1 a_{*} s_2 a_{*} s_3 p_1 s_4$ for some set partitions $s_1, s_2, s_3, $ and $s_4$. In addition, let $S = \{s_1, s_2, s_3, s_4\}$. Now, for any $s \in S$, if some $a \in s$ satisfies $|\mathrm{ind}_{s}(a)| = 2$, then $\mathrm{nc}(s)$ is clumped in $\phi_{aba}(p)$ by \Cref{lem1}. But this negates \Cref{cor2} for $i=1$. Therefore, all $a \in s$ for each $s \in S$ must satisfy $|\mathrm{ind}_{s}(a)| = 1$. Furthermore, if some $a \in s_1$ also satisfies $a \in s_2$ (resp. $s_3$), then $\mathrm{nc}(s_1 a_{*} s_2 a_{*} s_3)$ is clumped in $\phi_{aba}(p)$. But if so, then $C(\phi_{aba}(p)) > 1$, which negates \Cref{cor2} for $i=1$; thus, $s_1 \cap s_2 = \emptyset$ (resp, $s_1 \cap s_3 = \emptyset$). As a result, $\phi_{aba}(p) = \mathrm{r}(s_2) \mathrm{r}(s_3) {a_{*}}^2 \mathrm{r}(s_1) \mathrm{r}(s_4) p_1^2$ and so, $\mathrm{trunc}(\phi_{aba}(p)) = \mathrm{r}(s_2) \mathrm{r}(s_3) a_{*} \mathrm{r}(s_1) \mathrm{r}(s_4) p_1$.

    Next, by \Cref{trunc}, if $p$ is not sorted by $\phi_{aba}^{N(p)-1},$ then $\mathrm{r}(s_2)\mathrm{r}(s_3)a_{*}\mathrm{r}(s_1)\mathrm{r}(s_4)$ must not be sorted by $\phi_{aba}^{N(p)-2}$. Thus, by \Cref{main}, it must be that \[\mathrm{r}(s_2)\mathrm{r}(s_3)a_*\mathrm{r}(s_1)\mathrm{r}(s_4) = (\phi_{aba}(p)_1 \phi_{aba}(p)_2 \cdots \phi_{aba}(p)_{N(p)-1})^2.\]
    
    Now, because $a_{*} \in s_4$ by \Cref{theorem2lem1} and $s_1 \cap (s_2 \cup s_3) = \emptyset$, it must be that $|r(s_4)| \geq N(p)-1$. But because each $a \in s$ satisfies $|\mathrm{ind}_{a}(s_4)|=1$, it holds that $|r(s_4)| \leq N(p)-1$. Thus, $|r(s_4)| = N(p)-1$. As a result, $$\mathrm{r}(s_2) \mathrm{r}(s_3) a_{*} \mathrm{r}(s_1) = \mathrm{r}(s_4) = r(p_{\mathrm{ind}_{p}^{2}(p_1)+1}p_{\mathrm{ind}_{p}^{2}(p_1)+2} \ldots p_{|p|}).$$ Thus, fixing $|s_1|,$ $|s_2|,$ and $|s_3|$ such that $|s_1| +|s_2| + |s_3| = N(p)-1$ fixes $p$.
\end{proof}

Next, we count the number of inequivalent set partitions $p$ that satisfy $|p| = 2(N(p)-1) + |\mathrm{ind}_p(p_1)|$ and are not sorted after applying $\phi_{aba}^{N(p)-1}$. 

\begin{lemma} \label{theorem2lem4}
    The number of inequivalent set partitions $p$ that satisfy $|p| = 2(N(p)-1) + |\mathrm{ind}_p(p_1)|$ and are not sorted after applying $\phi_{aba}^{N(p)-1}$ is given by $${2N(p)+|\mathrm{ind}_p(p_1)|-3\choose |\mathrm{ind}_p(p_1)|-1}-|\mathrm{ind}_p(p_1)|{N(p)+|\mathrm{ind}_p(p_1)|-3 \choose |\mathrm{ind}_p(p_1)|-1}.$$ 

\end{lemma}
\begin{proof}
    By \Cref{cor2}, $C(\phi^0_{aba}(p)) = C(p) =0$. Thus, all $a (\neq p_1) \in p$ must satisfy $|\mathrm{ind}_p(a)|=2$. Let  $p=p_1 s_1 p_1 s_2 \cdots p_1 s_m$ where $|s_i| \geq 0$ for all $1\leq i\leq m$. Let $S = \{s_1, s_2, \ldots, s_m\}$. Now, if there exists some $s \in S$ and $a \in s$ such that $|\mathrm{ind}_{s}(a)| = 2$, then $\mathrm{nc}(s)$ is clumped in $\phi_{aba}(p)$. But if so, $C(\phi_{aba}(p)) > 1$, which negates \Cref{cor2} for $i=1$. Thus, each $a \in s$ must satisfy $|\mathrm{ind}_{s}(a)| = 1$. In particular, $|s| \leq N(p)-1$ for all $1 \leq i \leq m$. 

    Now, by \Cref{lem1}, $\phi_{aba}(p) = \mathrm{r}(s_1) \cdots \mathrm{r}(s_m) p_1^m$. Thus, $\mathrm{trunc}(\phi_{aba}(p)) = \mathrm{r}(s_1) \cdots \mathrm{r}(s_m) p_1$. Next, by \Cref{trunc}, if $p$ is not sorted by $\phi_{aba}^{N(p)-1},$ then $\mathrm{r}(s_1) \cdots \mathrm{r}(s_m)$ must not sorted by $\phi_{aba}^{N(p)-2}$. Thus, by \Cref{main}, it must be that $$\mathrm{r}(s_1) \cdots \mathrm{r}(s_m) = (\phi_{aba}(p)_1 \phi_{aba}(p)_2 \cdots \phi_{aba}(p)_{N(p)-1})^2.$$
    
    Therefore, fixing $|s_1|, |s_2|, \ldots, |s_m|$ such that $\sum_{i=1}^{m} |s_i| = 2N(p)-2$ and $|s_i| \leq N(p)-1$ for each $1 \leq i \leq m$ fixes $p$. 
\end{proof}

We end by using Lemmas \ref{char}, \ref{theorem2lem1}, \ref{theorem2lem2}, \ref{theorem2lem3}, and \ref{theorem2lem4} to prove \Cref{main2}.

\begin{proof}[Proof of \Cref{main2}]
   Let $a_{*}$ as in \Cref{char}. \Cref{theorem2lem1} shows that \Cref{theorem2lem2} and \Cref{theorem2lem3} count all of the set partitions that satisfy the statement of \Cref{main2} and $a_{*} \ne p_1$. In addition, setting $|\mathrm{ind}_{p}(p_1)|=3$ in \Cref{theorem2lem4} counts all set partitions that satisfy the statement of \Cref{main2} and $a_{*} = p_1$.
\end{proof}




\nocite{*}
\printbibliography

\end{document}